\documentclass[twoside,a4paper]{article}

\usepackage{amsmath,amsthm,amssymb}

\newcommand{\wap}{\operatorname{WAP}}
\newcommand{\ip}[2]{{\langle {#1} , {#2} \rangle}}
\newcommand{\aone}{\Box}
\newcommand{\atwo}{\Diamond}
\newcommand{\proten}{{\widehat{\otimes}}}
\newcommand{\mc}[1]{\mathcal{#1}}

\theoremstyle{plain}
\newtheorem{proposition}{Proposition}[section]
\newtheorem{theorem}[proposition]{Theorem}

\theoremstyle{definition}
\newtheorem{definition}[proposition]{Definition}

\theoremstyle{remark}

\begin{document}

\large
\title{\textsc{Weakly almost periodic functionals, representations, and operator spaces}}
\author{Matthew Daws}
\maketitle

\begin{abstract}
A theorem of Davis, Figiel, Johnson and Pe\l czy\'nski tells us that
weakly-compact operators between Banach spaces factor through reflexive
Banach spaces.  The machinery underlying this result is that of the real
interpolation method, which has been adapted to the category of operator
spaces by Xu, showing the this factorisation result also holds for
completely bounded weakly-compact maps.  In this note, we show that
Xu's ideas can be adapted to give an intrinsic characterisation of
when a completely contractive Banach algebra arises as a closed subalgebra
of the algebra of completely bounded operators on a reflexive operator
space.  This result was shown by Young for Banach algebras, and our
characterisation is a direct analogue of Young's, involving weakly
almost periodic functionals.
\end{abstract}



\section{Introduction}

In \cite{DFJP}, Davis, Figiel, Johnson and Pe\l czy\'nski showed that
weakly-compact operators and operators which factor through a reflexive
Banach space are the same class.  To do this, they used the real interpolation
space method, although this was not made explicit.  In \cite{Young} this result
was used to give a link between Banach algebras which arise as
closed subalgebra of the algebra of operators on a reflexive Banach space,
and Banach algebras which admit sufficiently many so called weakly almost
periodic functionals.  In \cite{Kai}, Kaiser noticed that this was
really a result about interpolation of Banach modules.  Recently, in
\cite{Daws}, we argued that such results can really be thought of as
results about representing dual Banach algebras.

The use of operator spaces has attracted a lot of attention in studying
non-self-adjoint operator algebras (see the monograph \cite{BLM} for
example).  Less well studied are
the much wider class of completely contractive (or quantum) Banach algebras.
The most common example is the Fourier algebra $A(G)$ for a locally
compact group $G$.  As $A(G)$ is the predual of the group von Neumann
algebra $VN(G)$, it carries a natural operator space structure.  It seems
that when $A(G)$ is considered as an operator space, properties of $A(G)$
better reflect properties of $G$ (see Ruan's original paper \cite{Ruan},
or the survey paper \cite{Runde2}).  This current paper was
motivated by concrete questions to do with $A(G)$, a matter we consider
further at the end of this paper.

I have learnt that F.\ Uygul has independently found the main result of
this paper (Theorem~\ref{Main_Thm} below).  See the paper \cite{Uygul}.
Our approaches are rather similar; this paper proves a factorisation result
for module maps, while Uygul takes an approach closer to the presentation of
\cite{Daws}.  However, the module result follows easily from Uygul's result.
We both essentially adapt ideas of Xu contained in \cite{Xu}.

As we suspect that this paper will mostly be read by people familiar
with Banach algebras, we shall follow the notation of \cite{Dales}, in
particular writing $E'$ for the dual of a Banach (or operator) space $E$,
and $T'$ for the adjoint of a linear map $T$.  We write $\kappa_E:E
\rightarrow E''$ for the canonical map from a Banach (or operator) space
to its dual.  With these exceptions,
we follow \cite{ER} for notation to do with operator spaces (see also
\cite{Pisier} for basic details about operator spaces).
A further exception is that when $T$ is an operator between operator
spaces $E$ and $F$, we always use brackets, and write $(T)_n$ for
the amplification map $\mathbb M_n(E) \rightarrow \mathbb M_n(F)$.
We shall also frequently abuse notation, and not notationally
distinguish between the norm on $\mathbb M_n(E)$ and that on
$\mathbb M_m(E)$.

\section{Completely contractive Banach algebras}

Let $\mc A$ be a Banach algebra which is also an operator space.
Let $\Delta:\mc A\otimes\mc A\rightarrow\mc A; a\otimes b\mapsto
ab$ be the multiplication map.  We say that $\mc A$ is a
\emph{completely contractive (CC) Banach algebra}
when $\Delta$ extends to a completely contractive map
$\mc A\proten\mc A\rightarrow \mc A$. Here, of course, $\proten$
denotes the operator space projective tensor product
(see \cite[Chapter~7]{ER}).  By using the identification
$\mc{CB}(\mc A\proten\mc A,\mc A) = \mc{CB}(\mc A,\mc{CB}(\mc A))$, we
see that $\mc A$ is a CC Banach algebra if and only
if the left-regular representation of $\mc A$ on itself maps into
$\mc{CB}(\mc A)$ and is a completely contractive homomorphism.

Similarly, let $\mc A$ be a CC Banach algebra, and let $E$ be a left $\mc A$-module.
Following \cite{Dales}, all our modules shall be contractive (which
can be arranged by a suitable renorming).  We say that $E$ is a
\emph{completely contractive (CC) left $\mc A$-module}
when $E$ is an operator space and the module map
$\mc A\otimes E\rightarrow E; a\otimes x\mapsto a\cdot x$ extends to
a complete contraction $\mc A\proten E\rightarrow E$.  As above,
this is equivalent to the induced homomorphism $\theta:\mc A\rightarrow
\mc B(E)$ actually mapping into $\mc{CB}(E)$ and being a complete contraction.
Obvious definitions apply to right $\mc A$-modules and $\mc A$-bimodules.

We turn $\mc A'$ into an $\mc A$-bimodule in the obvious way
\[ \ip{a\cdot \mu}{b} = \ip{\mu}{ba} = \ip{\mu\cdot b}{a}
\qquad (a,b\in\mc A, \mu\in\mc A'). \]
As $(\mc A\proten\mc A)' = \mc{CB}(\mc A,\mc A')$, it is easy to
see that when $\mc A$ is a CC Banach algebra, we have that $\mc A'$ is a CC
$\mc A$-bimodule, and also for $\mc A''$ and so forth.

The following is the operator space version of a notion first
formally defined by Runde in \cite{Runde1} (although it had been
studied before).

\begin{definition}
Let $E$ be an operator space, and suppose that $\mc A = E'$ is a CC Banach algebra.
When the product on $\mc A$ is separately weak$^*$-continuous, we
say that $\mc A$ is a \emph{completely contractive (CC) dual Banach algebra}.
\end{definition}

It is simple to check that the product is separately weak$^*$-continuous
if and only if $\kappa_E(E)$ becomes a sub-$\mc A$-bimodule of $\mc A'=E''$.
Notice that the operator space structure plays little role here,
essentially because duality works so well, \cite[Section~3.2]{ER}.
The equivalent notion, for algebras of operators on a Hilbert space,
has been widely studied (see \cite{LM1}, for example, where an operator version
of our result, Theorem~\ref{Main_Thm} below, is proved).

We showed in \cite{Daws} that the class of dual Banach algebras coincides
with the class of weak$^*$-closed subalgebras of $\mc B(E)$, where $E$ is a
reflexive Banach space.  This result follows quite easily from the
work of Young and Kasier, or at least their methods.  In this paper
we shall prove an analogous result for CC dual Banach algebra.  Notice that when
$E$ is a reflexive operator space, we have that $(E \proten E')' =
\mc{CB}(E)$, so that $\mc{CB}(E)$ is a dual space.  It is simple to
show that $\mc{CB}(E)$ is a CC Banach algebra, and that $E\proten E'$
is a submodule of $\mc{CB}(E)'$.  Thus $\mc{CB}(E)$ is a CC dual
Banach algebra.

\section{Interpolation spaces}

It seems that we cannot escape explaining a little about interpolation
spaces, given the important role played by the complex interpolation
method in the theory of operator spaces.

Let $E$ and $F$ be Banach spaces which are embedded continuously into
some Hausdorff topological vector space $X$ (in applications, we shall
typically have an injection of $E$ into $F$ allowing us to take $X=F$).
We define norms on the subspaces $E\cap F$ and $E+F$ of $X$ by
\[ \|x\|_{E\cap F} = \max\big( \|x\|_E , \|y\|_F \big) , \quad
\|w\|_{E+F} = \inf\big\{ \|x\|_E + \|y\|_F : w=x+y \big\}. \]
Obviously we can at this point replace $X$ by $E+F$ if we wish.  We
say that $(E,F)$ is a compatible couple.  Loosely, an interpolation
space is a Banach space intermediate to $E+F$ and $E\cap F$, and such
that certain mapping properties hold.  See \cite{BB} or \cite{BL}
for further details.

These ideas can be adapted, with minor tweaking, to the setting of
operator spaces, see \cite{Xu} or \cite[Section~2.7]{Pisier}.  The complex
interpolation method gives a Banach space $(E,F)_\theta$ where
$\theta$ is a parameter between $0$ and $1$.  When $E$ and $F$ are
operator spaces, $(E,F)_\theta$ becomes an operator space by setting
\[ \mathbb M_n\big( (E,F)_\theta \big) = \big( \mathbb M_n(E),
\mathbb M_n(F) \big)_\theta \qquad (n\geq1). \]
Let $(E_1,F_1)$ be another compatible couple, and let $T:E+F
\rightarrow E_1+F_1$ be a linear map such that $T$ maps $E$ into
$E_1$ in a (completely) bounded fashion, and the same for $F$ and
$F_1$.  Then $T:(E,F)_\theta\rightarrow(E_1,F_1)_\theta$ and $\|T\|
\leq \|T\|^{1-\theta}_{E\rightarrow E_1} \|T\|^\theta_{F\rightarrow F_1}$.
Let $\mu$ be a measure space and consider the Banach spaces $L_p(\mu)$.
Then $(L_1(\mu),L_\infty(\mu))$ is a compatible couple, and, isometrically,
\[ L_p(\mu) = \big(L_1(\mu),L_\infty(\mu)\big)_\theta
\qquad (1/p=1-\theta). \]
We give $L_\infty(\mu)$ the operator space structure it has as a
C$^*$-algebra, and give $L_1(\mu)$ the operator space structure it
gets from embedding $L_1(\mu)$ into the dual of $L_\infty(\mu)$.
Then the above identity allows us to define an operator space
structure on $L_p(\mu)$.  In particular, $L_2(\mu)$ gives an example
of Pisier's self-dual operator Hilbert space, see \cite[Section~7]{Pisier}.

Vector-valued versions of the above shall be important for us.
Let $(E_i)_{i\in I}$ be a family of operator spaces.  We let $\ell_\infty(E_i)$
be the usual direct sum of operator spaces, so that
\[ \ell_\infty(E_i) = \big\{ (x_i)_{i\in I} : x_i\in E_i \ (i\in I),
\|(x_i)\| := \sup_i \|x_i\| \big\}, \]
as a Banach space.  We define $\mathbb M_n(\ell_\infty(E_i)) =
\ell_\infty( \mathbb M_n(E_i) )$ for $n\geq 1$.  Similarly, $\ell_1(E_i)$
can be given an operator space structure by embedding it in $\ell_\infty(E_i')'$.
Alternatively, as noted in \cite[Section~2.6]{Pisier}, $\ell_1(E_i)$
is characterised by the universal property that whenever $E$ is an
operator space and $T_i:E_i\rightarrow E$ is a complete contraction,
for each $i\in I$, then the map $T:\ell_1(E_i)\rightarrow E$, given by
$T(x_i) = \sum_{i\in I} T_i(x_i)$, is a complete contraction.
When $I=\{1,2\}$, we write $E_1\oplus_1 E_2$ for $\ell_1(E_i)$, and
so forth.  Finally, we define
\[ \ell_p(E_i) = \big( \ell_1(E_i), \ell_\infty(E_i) \big)_\theta
\qquad (1/p=1-\theta). \]
We shall mainly use $\ell_2(E_i)$.  We note that when $E_i=\mathbb C$
for each $i$, we have that $\ell_2(E_i) = OH(I)$, Pisier's operator
Hilbert space, and that there are many characterisations of the
operator space structure on $OH(I)$.  We are not aware of characterisations
of $\ell_2(E_i)$ which do not use complex interpolation, however.

The real interpolation method is more complicated to explain, and
significantly harder to adapt to the operator space setting.  Fortunately,
Xu has done the hard work for us in \cite{Xu}.  We shall sketch a simple
case of the constructions Xu considers, following the approach of
Palmer in \cite[Section~1.7.8]{PalBook} and also as used by us in \cite{Daws}.

\section{Factoring module maps}

Firstly we shall consider the Banach space case.  Let $\mc A$ be
a Banach algebra, let $E$ and $F$ be left $\mc A$-modules, and let
$T:E\rightarrow F$ be an $\mc A$-module homomorphism.  That is,
$T(a\cdot x) = a\cdot T(x)$ for $a\in\mc A$ and $x\in E$.
For $n\in\mathbb N$, define a new norm $\|\cdot\|_n$ on $F$ by
\[ \|x\|_n = \inf\big\{ 2^{-n/2} \|T\|\|y\| + \|x-T(y)\| : y\in E \big\}
\qquad (x\in F). \]
We may check that
\[ \|x\|_n \leq \|x\| \leq 2^{n/2}\|x\|_n \qquad (x\in F), \]
so that $\|\cdot\|_n$ is an equivalent norm on $F$.  We let
\[ G = \Big\{ x\in F : \|x\|_G := \Big(\sum_{n=1}^\infty \|x\|_n^2\Big)^{1/2}
< \infty \Big\}, \]
so that for $y\in E$, as $\|T(y)\|_n \leq 2^{-n/2} \|T\| \|y\|$, we see that
$\|T(y)\|_G \leq \|T\| \|y\|$, from which it follows that $T(E)\subseteq G$.
Let $\iota:G\rightarrow F$ be the
inclusion map, which is norm-decreasing, and let $R$ be the map $T$,
treated as map from $E$ to $G$, so that $\|R\|\leq\|T\|$.  Hence $T$ factors
through the normed space $G$, as $\iota \circ R = T$.  Notice that $R$ has
dense range, and that $\iota$ is injective, so that $R$ has the same kernel
as $T$.

Let $E_0 = E / \ker T$, so that $T$ becomes an injection $E_0 \rightarrow
F$, and so we can regard $(E_0,F)$ as a compatible couple.  Then $G$
is a member of the equivalence class $(E_0,F)_{1/2,2}$, this being a
real interpolation space.  From standard results (see, for example,
\cite[Section~2.3, Proposition~1]{BB}) it follows that $G$ is reflexive
if and only if $T$ is weakly-compact.  In fact, this is not too hard
to prove directly, which we leave as an exercise for the reader.

Notice that for $x\in F$ and $a\in\mc A$, we see that $\|a\cdot x\|_n
\leq \|a\| \|x\|_n$ for $n\geq1$.  Hence $G$ becomes a left $\mc A$-module,
and $\iota$ and $R$ become $\mc A$-module homomorphisms.

We wish to carry out a similar construction for operator spaces.
A first step is to consider a different way of expressing $G$.  For
a Banach space $E$ and $t>0$, let $tE$ be the same space with the norm
multiplied by $t$ (and similarly for an operator space $E$, where we
set $\mathbb M_n(tE) = t\mathbb M_n(E)$).  Set
\[ F_n = 2^{-n/2}\|T\| E \oplus_1 F \qquad (n\geq 1), \]
let $Y = \ell_2(F_n)$, and let
\[ Y_0 = \{ (x_n,y_n)\in Y : T(x_n)+y_n=0 \ (n\geq 1) \}\subseteq Y. \]
It is easy to see that $Y_0$ is a closed subspace of $Y$, so we may
form the quotient space $Y/Y_0$.  Finally, define
\[ Y_1 = \{ (x_n,y_n)+Y_0 \in Y/Y_0 : (T(x_n)+y_n) \text{ constant} \}, \]
so that $Y_1$ is a closed subspace of $Y/Y_0$.

Define a map $\alpha:Y_1 \rightarrow G$ by
\[ \alpha( (x_n,y_n)+Y_0 ) = T(x_1)+y_1 \qquad ((x_n,y_n)+Y_0\in Y/Y_0). \]
It is easy to see that $\alpha$ is actually well-defined on the whole
of $Y/Y_0$, does map into $G$, and is an injection when restricted to $Y_1$.
Let $y\in G$, so for some sequence $(x_n)$ in $E$, we have that
\[ \sum_{n=1}^\infty \Big( 2^{-n/2}\|T\|\|x_n\| + \|y-T(x_n)\| \Big)^2
<\infty. \]
Let $y_n = y-T(x_n)$ for each $n$, so that $(x_n,y_n)\in Y$, and by
definition, $\alpha( (x_n,y_n)+Y_0 )=y$.  Hence $\alpha$ is a bijection.
A similar calculation shows that $\alpha$ is actually an isometry.

Now let $\mc A$ be a CC Banach algebra, let $E$ and $F$ be CC left $\mc A$-modules,
and let $T:E\rightarrow F$ be a completely bounded $\mc A$-module
homomorphism.  We note that, using the complex interpolation method
described above, we can give each $F_n$ and $Y$ natural operator space
structures.  Hence $Y_0, Y/Y_0$ and $Y_1$ all gain operator space
structures, and we can hence use $\alpha$ to induce an operator space
structure on $G$.  It follows from Xu's work that $\iota$ and $R$
are completely bounded.  Furthermore, as $G$ is a left $\mc A$-module,
we have a homomorphism $\theta_G:\mc A\rightarrow\mc B(G)$.  Xu's work
shows that $\theta_G$ actually maps into $\mc{CB}(G)$ and that $\theta_G$
is norm-decreasing.  However, we wish to show that $\theta_G$ is a
complete contraction, in order to show that $G$ is a CC left $\mc A$-module.
This follows from Xu's methods, but not seemingly directly from the results
of \cite{Xu}.

\begin{proposition}
With notation as above, there exists an absolute constant $K>0$
such that $\|R\|_{cb}\leq K\|T\|_{cb}$ and $\|\iota\|_{cb}\leq K$.
\end{proposition}
\begin{proof}
We shall only sketch this, as it follow from \cite[Theorem~2.2]{Xu},
and the remark thereafter.  By the universal property of $\oplus_1$, the map
\[ \iota_n: F_n = 2^{-n/2}\|T\|_{cb} E \oplus_1 F \rightarrow F;\
(x,y) \mapsto 2^{-n/2}(T(x)+y) 
\qquad (n\geq1, x\in E, y\in F) \]
is a complete contraction.  Now consider the maps, for $1\leq t\leq\infty$,
\[ \beta_t: \ell_t(F_n) \rightarrow F; \
(w_n) \mapsto \sum_{n=1}^\infty 2^{-n/2}\iota_n(w_n)
= \sum_{n=1}^\infty 2^{-n}(T(x_n)+y_n)
\qquad ( w_n=(x_n,y_n)\in F_n ). \]
By the universal property, as each
$\iota_n$ is a complete contraction, we have that $\|\beta_1\|_{cb}
\leq 2^{-1/2}$.  Let $(w_n)\in\mathbb M_k(\ell_\infty(F_n))$, so that
\[ \|\beta_\infty(w_n)\| \leq \sum_{n=1}^\infty
2^{-n/2}\|\iota_n(w_n)\| \leq \sum_{n=1}^\infty
2^{-n/2}\|w_n\|_{\mathbb M_k(F_n)} \leq 
\frac{\sqrt2}{2-\sqrt2} \|(w_n)\|, \]
so that $\|\beta_\infty\|_{cb}\leq2^{1/2}(2-2^{1/2})^{-1}$.
By complex interpolation,
\[ \|\beta_2\|_{cb} \leq \big(2^{-1/2}\big)^{1/2}
\big( 2^{1/2}(2-2^{1/2})^{-1} \big)^{1/2}
= \sqrt{ \frac{1}{2-\sqrt2} } = K. \]
Notice that $\beta_2$ vanishes on $Y_0$, and so $\beta_2$ drops to
a well-defined operator from $Y/Y_0$ to $F$, and hence by restriction
to a map $Y_1\rightarrow F$.  Fairly obviously, $\beta_2\alpha^{-1}
= \iota$, showing that $\iota$ is completely bounded.
The argument for $R$ follows similarly.
\end{proof}

Notice that
\[ \|\beta_2(w_n)\| \leq \sum_{n=1}^\infty 2^{-n/2} \|\iota_n(w_n)\|
\leq \Big(\sum_{n=1}^\infty \|\iota_n(w_n)\|^2 \Big)^{1/2}, \]
so we see that $\|\beta_2\|\leq 1$, so that the complex interpolation
estimate of $\|\beta_2\|$ is not optimal.  It seems possible that if
we had a more concrete description of the operator space structure on
$\ell_2(F_n)$, then we could show that $\|R\|_{cb} \|\iota\|_{cb}
\leq \|T\|_{cb}$.  However, the following is a suitable work around.

\begin{proposition}
Let $E, F$ and $G$ be operator spaces, and let $T\in\mc{CB}(E,F),
R\in\mc{CB}(E,G)$ and $\iota\in\mc{CB}(G,F)$ be such that $T=\iota R$.
We may give $G$ an equivalent operator space structure for which
$\iota$ is a complete contraction, and $\|R\|_{cb}=\|T\|_{cb}$.
\end{proposition}
\begin{proof}
By replacing $G$ by $tG$ for some $t>0$, we may suppose that $\iota$ is
already a complete contraction.  Define a map $\phi:\|T\|_{cb}E \oplus_1 G
\rightarrow G$ by $\phi(x,y) = R(x)+y$ for $x\in E$ and $y\in G$.
Clearly $\phi$ is surjective, so we may identify $G$ with
$\|T\|_{cb}E \oplus_1 G / \ker\phi$.  Use this to induce a new operator
space structure on $G$, say with norm $\|\cdot\|_0$, giving $G_0$.
Let $w\in\mathbb M_n(G)$, so that clearly $\|w\|_0 \leq \|w\|$.
Conversely, let $x\in\mathbb M_n(E)$ and $y\in\mathbb M_n(G)$ be
such that $w=(R)_n(x)+y$.  The map $E\rightarrow G; x \mapsto
\|T\|_{cb} \|R\|_{cb}^{-1} R(x)$ is a complete contraction, so by the
universal property of $\oplus_1$,
\[ \|w\| = \big\|(R)_n(x) + y\big\|
\leq \big\| \big( \|T\|_{cb}^{-1}\|R\|_{cb}x , y \big) \big\|. \]
As $\|T\|_{cb} \leq \|R\|_{cb}$, the map $\|T\|_{cb}E\oplus_1 G
\rightarrow \|T\|_{cb}E\oplus_1 G; (x,y) \mapsto (x,\|T\|_{cb}\|R\|_{cb}^{-1})$
is a complete contraction, so that
\[ \|w\| \leq \|T\|_{cb}^{-1} \|R\|_{cb} \|(x,y)\|, \]
from which we conclude that $\|w\|\leq\|T\|_{cb}^{-1} \|R\|_{cb}\|w\|_0$,
showing that $\|\cdot\|_0$ is an equivalent operator space structure on $G$.
Clearly we have that $\|R:E\rightarrow G_0\|_{cb} \leq \|T\|_{cb}$.

As $\phi$ is defined to be a complete quotient map, $\phi':G_0'\rightarrow
\|T\|_{cb}^{-1}E' \oplus_\infty G'$ is a complete isometry.  We may check that
\[ \phi'(\mu) = \big( R'(\mu), \mu \big)
\qquad (\mu\in G_0'). \]
Let $\lambda\in \mathbb M_n(F')$, so that
\begin{align*}
\|(\iota')_n(\lambda)\|_{G_0'} &= \big\| ( (R')_n(\iota')_n(\lambda) ,
   (\iota')_n(\lambda) ) \big\|_{\|T\|_{cb}^{-1}E' \oplus_\infty G'} \\
&= \max\big( \|T\|_{cb}^{-1} \|(T')_n(\lambda)\| , \|(\iota')_n(\lambda)\| \big)
\leq \|\lambda\|.
\end{align*}
Hence $\| \iota:G_0\rightarrow F \|_{cb}\leq 1$, as required.
\end{proof}

Recall the homomorphism $\theta_G:\mc A\rightarrow\mc B(G)$, which
actually maps into $\mc{CB}(G)$ by Xu's work.

\begin{proposition}\label{mod_okay}
With notation as above, $\theta_G$ is a complete contraction.
\end{proposition}
\begin{proof}
We first show that our claim holds for the original definition of
$G$, and then check that the above renorming procedure does not chance
our conclusions.
Fix $n\geq 1$ and $a=(a_{kl})\in\mathbb M_n(\mc A)$ with $\|a\|\leq 1$.
Then $A:=(\theta_G)_n(a)\in\mathbb M_n(\mc{CB}(G)) =
\mc{CB}(G,\mathbb M_n(G))$.  Let $m\geq 1$ and let $w\in\mathbb M_m(G)$,
so that $(A)_m(w)\in\mathbb M_m(\mathbb M_n(G))
= \mathbb M_{m\times n}(G)$.  We wish to show that $\|(A)_m(w)\| \leq \|w\|$,
which would both demonstrate Xu's result that $\theta_G$ maps into
$\mc{CB}(G)$, and would show that $\theta_G$ is a complete contraction.

Recall the definitions of $(F_k)_{k\geq 1}, Y, Y_0, Y_1$ and $\alpha$
from above.
Clearly, for $k\geq 1$,
\[ (\theta_E)_n(a)\in\mc{CB}\big( 2^{-k/2}\|T\|_{cb}E,
\mathbb M_n(2^{-k/2}\|T\|_{cb}E) \big) \]
with norm at most $\|a\|\leq1$, as $\theta_E$ is a complete contraction.
By the universal property of $\oplus_1$, it is clear that
\[ A_k:= (\theta_E)_n(a) \oplus (\theta_F)_n(a) : F_k \rightarrow
\mathbb M_n(F_k) \]
is a complete contraction, for each $k\geq 1$.  Clearly the diagonal
map $(\oplus_k A_k):\ell_\infty(F_k) \rightarrow \mathbb M_n(\ell_\infty(F_k))
= \ell_\infty(\mathbb M_n(F_k))$ is a complete contraction.  Again, by
the universal property of $\oplus_1$, $(\oplus_k A_k):\ell_1(F_k)
\rightarrow \mathbb M_n(\ell_1(F_k))$ is also complete contraction.
Hence, by complex interpolation,
\[ (\oplus_k A_k): Y=\ell_2(F_k) \rightarrow \mathbb M_n(Y) \]
is a complete contraction.  It is clear that as $T$ is an $\mc A$-module
homomorphism, $(\oplus_k A_k)$ leaves $Y_0$ invariant, and hence
$(\oplus_k A_k)$ drops to a complete contraction on $Y/Y_0$.  Similarly,
$(\oplus_k A_k)$ restricts to a complete contraction on $Y_1$.

It is a simple check that $(\alpha)_{m\times n} (\oplus_k A_k)
(\alpha^{-1})_m = (A)_m$, and so we conclude that $(A)_m$ is a
contraction, as required.

Now consider the renorming.  Again let $a\in\mathbb M_n(\mc A)$ and
$x\in\mathbb M_m(G)$, so that $(A)_m(x)\in\mathbb M_{m\times n}(G)$.
Suppose that $\|x\|_0<1$, so that $x=\phi(y,z) = (R)_m(y)+z$ for some
$y\in\mathbb M_m(E)$ and $z\in\mathbb M_m(G)$ with
$\|(y,z)\|_{\|T\|_{cb}E\oplus_1 G} < 1$.  Let
\[ w = \big( (A)_m(y) , (A)_m(z) \big) \in
\mathbb M_m( \|T\|_{cb}E \oplus_1 G ), \]
so that
\[ \|w\|_{\|T\|_{cb}E\oplus_1 G} \leq
\|A\|_{cb}\|y,z\|_{\|T\|_{cb}E\oplus_1 G} < \|a\|. \]
We then observe that
\[ (A)_m(x) = ((\theta_G)_n(a))_m(R)_m(y) + A(z)
= (R)_{m\times n} ((\theta_E)_n(a))_m(y) + A(z)
= (\phi)_n(w), \]
so that $\|(A)_m(x)\| \leq \|\phi\|_{cb} \|w\| < \|a\|$, showing
that $\|A:G_0\rightarrow G_0\|_{cb}\leq\|a\|$ as required.
\end{proof}

We have hence shown the following factorisation result.

\begin{theorem}\label{mod_fac_thm}
Let $\mc A$ be a CC Banach algebra, let $E$ and $F$ be CC left $\mc A$-modules,
and let $T:E\rightarrow F$ be a completely bounded $\mc A$-module
homomorphism.  The following are equivalent:
\begin{enumerate}
\item $T$ is weakly-compact;
\item there exists a reflexive CC left $\mc A$-module $G$, an injective
complete contraction $\iota:G\rightarrow F$, and a completely bounded
map $R:E\rightarrow G$ with $\|R\|_{cb} = \|T\|_{cb}$, $T = \iota R$,
and such that $R$ and $\iota$ are $\mc A$-module homomorphisms.
\end{enumerate}
\end{theorem}

\section{Representing CC dual Banach algebras}

Theorem~\ref{mod_fac_thm} is central to proving Young's representation
theorem (a fact explicitly noticed by Kaiser in \cite{Kai}).  The
situation for CC Banach algebras is more complicated,
because we need to take account of the matrix structures $\mathbb M_n(\mc A)$.
Recall that while $\mathbb M_n(\mc A)$ is obviously an algebra, the
product need not be uniformly bounded in $n$.  Indeed, this is equivalent
to $\mc A$ being completely isomorphic to a subalgebra of $\mc B(H)$ for
a Hilbert space $H$ (see \cite[Chapter~17]{ER} and \cite{LM1} for example).

\begin{theorem}\label{Main_Thm}
Let $\mc A$ be a CC dual Banach algebra with predual $\mc A_*$.
Then there exists a weak$^*$-weak$^*$-continuous completely
isometric homomorphism $\theta:\mc A\rightarrow\mc{CB}(E)$ for
some reflexive operator space $E$.
\end{theorem}
\begin{proof}
We may suppose that $\mc A$ is unital.  If not, we replace $\mc A$ by
$\mc A\oplus_1 \mathbb C$, the unitisation of $\mc A$, and we may check
that properties of $\oplus_1$ imply that $\mc A\oplus_1 \mathbb C$ is
a CC dual Banach algebra, with predual $\mc A_* \oplus_\infty \mathbb C$.
Denote the unit of $\mc A$ by $e_{\mc A}$.

Let $n\geq 1$, we have that $\mathbb M_n(\mc A_*)$ is a CC left
$\mc A$-module.  Let $\mu=(\mu_{ij})\in\mathbb M_n(\mc A_*)$ with $\|\mu\|=1$,
and let $T:\mc A \rightarrow \mathbb M_n(\mc A_*)$ be the map
$T(a) = a\cdot\mu$, so that $T$ is completely contractive.  Recall that
Gantmacher's Theorem tells us that $T$ is weakly-compact if and only if
$T''(\mc A'') \subseteq \mathbb M_n(\mc A_*)$.  Let $b=(b_{ij}) \in \mathbb
M_n(\mc A_*)' = \mathbb T_n(\mc A)$ (see \cite[Proposition~7.1.6]{ER}).  Then
\[ \ip{T'(b)}{a} = \sum_{i,j=1}^n \ip{b_{ij}}{a\cdot\mu_{ij}}
= \sum_{i,j=1}^n \ip{\mu_{ij}\cdot b_{ij}}{a}, \]
so that $T'(b) = \sum_{i,j} \mu_{ij}\cdot b_{ij} \in\mc A_* \subseteq\mc A'$.
For $\Phi\in\mc A''$, let $c\in\mc A$ be such that $\ip{c}{\mu} =
\ip{\Phi}{\mu}$ for $\mu\in\mc A_*$.  Then we have that
\[ \ip{T''(\Phi)}{b} = \sum_{i,j=1}^n \ip{\Phi}{\mu_{ij}\cdot b_{ij}}
= \sum_{i,j=1}^n \ip{c}{\mu_{ij}\cdot b_{ij}}
= \sum_{i,j=1}^n \ip{b_{ij}}{c\cdot\mu_{ij}}. \]
Thus $T''(\Phi) = T(c)\in\mathbb M_n(\mc A_*)$, and we conclude
that $T$ is indeed weakly-compact.

Applying Theorem~\ref{mod_fac_thm}, we find a reflexive operator
space $G_\mu$ and complete contractions $\iota_\mu:G_\mu\rightarrow
\mathbb M_n(\mc A_*)$ and $R_\mu:\mc A\rightarrow
G_\mu$, such that $\iota_\mu R_\mu = T$, and with $\iota_\mu$ and $R_\mu$
being $\mc A$-module homomorphisms.  Let $\theta_\mu:\mc A\rightarrow
\mc{CB}(G_\mu)$ be induced by the module action.  Let $x_\mu=R_\mu(e_{\mc A})
\in G_\mu$, so that
\[ T(a) = \iota_\mu R_\mu(a\cdot e_{\mc A}) = 
\iota_mu \theta_\mu(a) R_\mu(e_{\mc A}) =
\iota_\mu \theta_\mu(a)(x_\mu)
\qquad (a\in\mc A). \]
Hence, for $a=(a_{kl})\in\mathbb M_n(\mc A)$,
\[ (T)_n(a) = \big( a_{kl}\cdot\mu \big)
= \big( \iota_\mu \theta_\mu(a_{kl})(x_\mu) \big)
= (\iota_\mu)_n \big( \theta_\mu(a_{kl})(x_\mu) \big)
= (\iota_\mu)_n (\theta_\mu)_n(a)(x_\mu), \]
noting that $(\theta_\mu)_n(a) \in \mc{CB}(G,\mathbb M_n(G))$.
Then note that
\begin{align*} \langle\ip{a}{\mu}\rangle
&= \big( \ip{a_{kl}}{\mu_{ij}} \big)
= \big( \ip{e_{\mc A}}{a_{kl}\cdot\mu_{ij}} \big)
= \big((T)_n(a)\big)(e_{\mc A})
= \big( (\iota_\mu)_n (\theta_\mu)_n(a)(x_\mu) \big) (e_{\mc A}),
\end{align*}
so that, as $\|x_\mu\| \leq 1$,
\[ \|(\theta_\mu)_n(a)\| \geq \| (\theta_\mu)_n(a)(x_\mu) \|
\geq \big\| (\iota_\mu)_n \big( (\theta_\mu)_n(a)(x_\mu) \big) (e_{\mc A}) \big\|
\geq \| \langle \ip{\mu}{a} \rangle \|. \]

Let $\Lambda = \{ \mu\in\mathbb M_n(\mc A_*) : n\geq 1, \|\mu\|=1 \}$,
and let $E = \ell_2(\{G_\mu : \mu\in\Lambda\})$.  Using the same argument as
in the proof of Proposition~\ref{mod_okay}, we see that the homomorphism
$\theta = \oplus_{\mu\in\Lambda} \theta_\mu:\mc A\rightarrow\mc B(E)$ maps
into $\mc{CB}(E)$, and is a complete contraction.  We can treat $x_\mu$
as a member of $E$ for each $\mu$, so we see that for $a\in\mathbb M_n(\mc A)$,
\begin{align*} \|a\| &=
\sup\{ \| \langle \ip{\mu}{a} \rangle \| : \mu\in\mathbb M_n(\mc A_*), \|mu\|=1 \} \\
&\leq \sup\{ \| (\theta_\mu)_n(a) \| : \mu\in\mathbb M_n(\mc A_*), \|\mu\|=1 \}
\leq \|(\theta)_n(a)\| \leq \|a\|, \end{align*}
and so we must have that $\theta$ is a complete isometry, as required.

Finally, we wish to show that $\theta$ is weak$^*$-weak$^*$-continuous.
That is, we wish to show that there exists a complete contraction
$\theta_*:E'\proten E\rightarrow\mc A_*$ such that $\theta_*'=\theta$.
As $\theta$ is a diagonal map on $E = \ell_2(G_\mu)$, it is sufficient
to check that each $\theta_\mu$ is weak$^*$-continuous.  Fix $\mu=(\mu_{ij})
\in\mathbb M_n(\mc A_*)$.  Recall that $T_\mu$ has dense range, and that
$\iota_\mu$ is injective.  We have that $\iota_\mu':\mathbb T_n(\mc A)
\rightarrow G_\mu'$ has dense range if and only if $\iota_\mu''$ is injective.
Now, as $G_\mu$ is reflexive, we identify $G_\mu$ with $G_\mu''$, and
we see that $\iota_\mu'' = \kappa_{\mathbb M_n(\mc A_*)}\iota_\mu$,
which is injective, showing that $\iota_\mu'$ does indeed have dense range.
Let $x\in G_\mu$ and $\lambda\in G_\mu'$.  By density, we may suppose that
$x = R_\mu(b)$ and that $\lambda=\iota_\mu'\kappa_{\mc A}(c)$ for some
$b\in\mc A$ and $c=(c_{ij})\in\mathbb T_n(\mc A)$.  Then
\begin{align*}
\ip{\theta_\mu'\kappa_{E_\mu'\proten E_\mu}(\lambda\otimes x)}{a}
&= \ip{\lambda}{\theta_\mu(a)(x)}
= \ip{\iota_\mu'\kappa_{\mc A}(c)}{\theta_\mu(a)R_\mu(b)}
= \ip{\iota_\mu R_\mu(ab)}{c} \\ &= \ip{ab\cdot\mu}{c}
= \sum_{i,j=1}^n \ip{ab\cdot\mu_{ij}}{c_{ij}}
= \sum_{i,j=1}^n \ip{a}{b\cdot\mu_{ij}\cdot c_{ij}},
\end{align*}
so that $\theta_\mu'\kappa_{E_\mu'\proten E_\mu}(\lambda\otimes x) \in \mc A_*$.
Hence $\theta_\mu$ is weak$^*$-continuous, as required.
\end{proof}

\subsection{Representing general CC Banach algebras}

Let $\mc A$ be a Banach algebra, let $E$ be a left $\mc A$-module,
and define
\[ \wap(E,\mc A) = \wap(E) = \{ x\in E : \mc A\rightarrow E;
a\mapsto a\cdot x \text{ is weakly-compact} \}, \]
and similarly for CC Banach algebras $\mc A$ and CC left $\mc A$-modules.
In particular, we say that $\wap(\mc A')\subseteq \mc A'$ is the
space of \emph{weakly almost periodic functionals}, the term coming
from abstract harmonic analysis.  We note that some authors write
$\wap(\mc A)$ for this space.

Let us quickly recall the Arens products (see \cite[Section~2]{Daws}
for further details about the following ideas, although we note that
most of this is folklore).  Let $\mc A$ be a Banach algebra.
We define bilinear actions $\mc A''\times
\mc A', \mc A'\times\mc A'' \rightarrow \mc A'$ by
\[ \ip{\Phi\cdot\mu}{a} = \ip{\Phi}{\mu\cdot a}, \quad
\ip{\mu\cdot\Phi}{a} = \ip{\Phi}{a\cdot\mu}
\qquad (\Phi\in\mc A'', \mu\in\mc A', a\in\mc A). \]
Then we define bilinear maps $\aone,\atwo:\mc A''\times\mc A''
\rightarrow\mc A''$ by
\[ \ip{\Phi\aone\Psi}{\mu} = \ip{\Phi}{\Psi\cdot\mu},\quad
\ip{\Phi\atwo\Psi}{\mu} = \ip{\Psi}{\mu\cdot\Phi}
\qquad (\Phi,\Psi\in\mc A'',\mu\in\mc A'). \]
It can be shown that $\aone$ and $\atwo$ are Banach algebra products on
$\mc A''$, called the \emph{first} and \emph{second Arens products},
and that $a\cdot\Phi = \kappa_{\mc A}(a)\aone\Phi =
\kappa_{\mc A}(a)\atwo\Phi$ for $a\in\mc A,\Phi\in\mc A''$, and
similarly on the right.  When $\aone=\atwo$, we say that $\mc A$ is
\emph{Arens regular}.  This is equivalent to $\wap(\mc A')=\mc A'$.
Indeed, more is true, as a shall see shortly.

We now sketch how to apply these ideas to a CC Banach algebra $\mc A$.
As before, this idea has been studied for operator algebras (see
\cite{Ruan1} for example) but we have not been able to find a good
source for CC Banach algebras; no doubt the following is known to
experts though.  Let $\Theta:\mc A\proten\mc A\rightarrow\mc A$ be
the completely contractive multiplication map, so that $\Theta':\mc A'
\rightarrow\mc{CB}(\mc A,\mc A')$ is also a complete contraction.
With the convention that
\[ \ip{T}{a\otimes b} = \ip{T(a)}{b}
\qquad (T\in(\mc A\proten\mc A)'=\mc{CB}(\mc A,\mc A'),
a\otimes b\in\mc A\proten\mc A), \]
we may check that $\Theta'(\mu)(a) = \mu\cdot a$ for $a\in\mc A$
and $\mu\in\mc A'$.  Define a map $\alpha:\mc A''\proten\mc A''
\rightarrow(\mc A\proten\mc A)''=\mc{CB}(\mc A,\mc A')'$ by
\[ \ip{\alpha(\Phi\otimes\Psi)}{T} = \ip{\Phi}{T'(\Psi)} 
= \ip{T''(\Phi)}{\Psi}
\qquad (\Phi,\Psi\in\mc A'', T\in\mc{CB}(\mc A,\mc A')), \]
so that as $(\mc A''\proten\mc A'')'=\mc{CB}(\mc A'',\mc A''')$
and the map $\mc{CB}(\mc A,\mc A')\rightarrow \mc{CB}(\mc A'',\mc A''');
T\mapsto T''$ is a complete isometry, we see that $\alpha$ is a
complete contraction.  Finally, we see that
\[ \ip{\Theta''\alpha(\Phi\otimes\Psi)}{\mu}
= \ip{\Phi}{\Theta'(\mu)'(\Psi)} = \ip{\Phi}{\Psi\cdot\mu}
= \ip{\Phi\aone\Psi}{\mu}
\qquad (\Phi,\Psi\in\mc A'', \mu\in\mc A'), \]
which shows that $(\mc A'',\aone)$ is a CC Banach algebra.
If we choose the other convention for identifying $\mc{CB}(\mc A,
\mc A')$ with the dual of $\mc A\proten\mc A$, we will find a proof
that $(\mc A'',\atwo)$ is a CC Banach algebra.

\begin{proposition}
Let $\mc A$ be a Banach algebra, and let $X\subseteq\mc A'$ be a closed
submodule.  The following are equivalent:
\begin{enumerate}
\item $X\subseteq\wap(\mc A')$;
\item the first (or, equivalently, second) Arens product drops to a
well-defined product on $X' = \mc A''/X^\perp$ turning $(X',X)$ into
a dual Banach algebra.
\end{enumerate}
\end{proposition}
\begin{proof}
See, for example, \cite[Proposition~2.4]{Daws} or compare
with \cite[Theorem~5.6]{Lau}.
\end{proof}

Loosely, we can say that $\mc A''/\wap(\mc A')^\perp$ is the largest
quotient of $\mc A''$ on which the Arens products agree.  The above will
clearly still hold for CC Banach algebras.  Combining these observations
with our factorisation theorem, we have the following.

\begin{theorem}\label{wap_thm}
Let $\mc A$ be a CC Banach algebra, and let $q:\mc A''
\rightarrow \wap(\mc A')'=\mc A''/\wap(\mc A')^\perp$ be the quotient
map.  The following are equivalent:
\begin{enumerate}
\item\label{cond_one} the map $q\kappa_{\mc A}:\mc A\rightarrow\wap(\mc A')'$ is
a complete isometry;
\item\label{cond_two} there is a homomorphism $\theta:\mc A\rightarrow \mc{CB}(E)$,
which is a complete isometry, for some reflexive Banach space $E$.
\end{enumerate}
We may replace the word ``isometry'' by ``isomorphism onto its
range'' above.  Furthermore, we may also replace the phrase ``a complete isometry''
by ``an injection'' above.
\end{theorem}
\begin{proof}
We shall show the isometric version; the isomorphic and injective
versions are similar.
If (\ref{cond_one}) holds, then $\mc A$ is completely isometric to
a subalgebra of some CC dual Banach algebra, and so (\ref{cond_two})
holds by Theorem~\ref{Main_Thm}.

Conversely, suppose there exists a complete
isometry $\theta:\mc A\rightarrow\mc{CB}(E)$.  For $n,m\geq 1$,
$x\in\mathbb M_m(E)$ and $\mu=(\mu_{ij})\in\mathbb M_{n\times m}(E')$ with
$\|x\|=\|\mu\|=1$, consider the map
\[ \Lambda_{x,\mu} = \Lambda:\mathbb M_n(\mc A)\rightarrow\mathbb M_{(n\times m)^2}; \
a=(a_{kl})\mapsto \langle\ip{\mu}{((\theta)_n(a))_m(x)}\rangle. \]
Then notice that for $a=(a_{kl})\in\mathbb M_n(\mc A)$,
\[ \Lambda(a) = \big( \ip{\mu_{ij}}{\theta(a_{kl})(x_{rs})} \big)
= \big( \ip{\theta'(\mu_{ij}\otimes x_{rs})}{a_{kl}} \big)
= \langle\ip{\big(\theta'(\mu_{ij}\otimes x_{rs})\big)}{a}\rangle
= \langle\ip{\lambda}{a}\rangle, \]
say, for some $\lambda\in\mathbb M_{n\times m\times m}(\mc A')$.
We claim that actually $\lambda\in\mathbb M_{n\times m\times m}(\wap(\mc A'))$,
which is equivalent to $\theta'(\mu_{ij}\otimes x_{rs}) \in \wap(\mc A')$
for each $i,j,r,s$.  The claim follows by the observation that the map
$\mc A\rightarrow\mc A'; a\mapsto a\cdot\theta'(\phi\otimes y)$ factors through
the reflexive Banach space $E$, for any $\phi\in E'$ and $y\in E$.
For $b\in\mathbb M_{n\times m\times m}(\mc A)$, we have that
\[ \langle\ip{\lambda}{b}\rangle
= \big(\ip{\theta'(\mu_{ij}\otimes x_{rs})}{b_{kl}}\big)
= \big( \ip{\mu_{ij}}{\theta(b_{kl})(x_{rs})} \big)
= \langle\ip{\mu}{((\theta)_{n\times m\times m}(b))_m(x)}\rangle, \]
We thus see that
\[ \|\lambda\| = \sup\{ \| \langle\ip{\lambda}{b}\rangle \| :
b\in\mathbb M_{n\times m\times m}(\mc A), \|b\|\leq 1 \}
\leq \|\mu\| \|x\| = 1. \]
We conclude that for $a\in\mathbb M_n(\mc A)$,
\begin{align*} \| (q\kappa_{\mc A})_n(a) \|
&= \sup\{ \|\langle\ip{\lambda}{a}\rangle\| : \lambda\in
   \mathbb M_{n\times m\times m}(\wap(\mc A')),
   \|\lambda\|\leq 1 \} \\
&\geq \sup\{ \|\Lambda_{x,\mu}(a)\| : \|x\|=\|\mu\|=1 \}
= \|(\theta)_n(a)\|_{cb} = \|a\|,
\end{align*}
while clearly $\| (q\kappa_{\mc A})_n(a) \| \leq \|a\|$, showing that
$q\kappa_{\mc A}$ is a complete isometry, as required.
\end{proof}

In the case of Banach algebras, the above is due to Young,
\cite{Young}; our proof is closer in nature to Kaiser's presentation
in \cite{Kai}.

Notice that $q\kappa_{\mc A}$ being an injection is independent of
the particular operator space structure on $\mc A$; this can be
restated by saying that a CC Banach algebra admits an injective
representation on a reflexive operator space if and only if the
underlying Banach algebra admits an injective representation on
a reflexive Banach space.  In contrast, it seems possible that
$q\kappa_{\mc A}$ might be an isomorphism, while $\mc A$ admits some
operator space structure turning it into a CC Banach algebra for which
$q\kappa_{\mc A}$ is not a \emph{complete} isomorphism.  We consider
this question for Fourier algebras below.

\section{Fourier algebras}

Let $G$ be a locally compact group, and consider the convolution algebra
$L_1(G)$.  The space $\wap(L_1(G)')$ is a classical object, which has
been widely studied (along with generalisations for semigroups, see
\cite{BJM}).  In particular, $\wap(L_1(G)')$ is a sub-C$^*$-algebra of
$L_\infty(G)$, say with (compact) character space $G^{\wap}$.  Then $G$
naturally embeds densely into $G^{\wap}$, and $G^{\wap}$ inherits a semigroup
structure (in fact induced by the Arens products, see \cite[Section~7]{Daws} for
example).  We can abstractly characterise $G^{\wap}$ as a certain semigroup
compactification of $G$.

Now consider instead the Fourier algebra $A(G)$.  When $G$ is abelian,
$A(G) = L_1(\hat G)$ where $\hat G$ is the dual group of $G$.  Hence
$\wap(A(G)')$ has an interpretation in terms of $\hat G$.  In particular,
$\wap(A(G)')$ is a sub-C$^*$-algebra of $VN(G)$, and certainly
$\wap(A(G)')$ is not all of $VN(G)$ (as $L_1(G)$ is only Arens regular
when $G$ is finite, indeed, see \cite{LL} which proves much more).
For non-abelian groups, $\wap(A(G)')$ was first studied in \cite{Gran}
and \cite{DR}.  Surprisingly, when $G$ is not abelian, it is, in general,
unknown if $\wap(A(G)')$ is a sub-C$^*$-algebra of $VN(G)$ (see, for
example, \cite{Hu} for recent work on this problem).
Indeed, it is not even known if $\wap(A(G)')=VN(G)$ can occur for infinite
$G$ (see \cite{Forrest1} for partial results).
However, when $G$ is amenable
and discrete, for example, it is known that $\wap(A(G)') = C^*_r(G)$, the
reduced C$^*$-algebra of $G$, as we would expect by analogy with the
abelian case.

In the context of Kac algebras, see \cite{ES}, or Locally Compact Quantum
Groups, see \cite{KV}, we view $L_1(G)$ and $A(G)$ as being dual to each
other, in some technical sense.  As noted in the introduction, it seems
to be necessary to use the operator space structure on $A(G)$ to fully
realise this idea.  Indeed, one could also argue that one should think
of $L_1(G)$ as an operator space, but as its dual is a commutative
C$^*$-algebra, $L_1(G)$ gets the max quantisation, and so in this case
we actually do not gain any new structure over viewing $L_1(G)$ as
simply a Banach space.

Actually, this is not quite true.  For example, let $\mc A$ be a
closed subalgebra of $\mc B(E)$ for some reflexive Banach space $E$.
It would seem to be a reasonable conjecture that if we give $\mc A$
the max quantisation (see \cite[Section~3.3]{ER}) then $\mc A$ becomes
a closed subalgebra of $\mc{CB}(F)$ for a suitable reflexive operator
space $F$.  We have, however, been unable to prove this\footnote{As
Uygul implicitly points out in \cite{Uygul}, this is actually rather simple
for dual Banach algebras.  Let $\mc A$ be a dual Banach algebra
with predual $\mc A_*$, and give $\mc A_*$ the min quantisation, so that
$\mc A$ gets the max quantisation, and is hence a CC (dual) Banach alegbra.
Hence $\mc A$ is weak$^*$-weak$^*$ isometric to a weak$^*$-closed subalgebra
of $\mc{CB}(E)$ for some reflexive operator space $E$.  Morally, we should
now be able to draw conclusions for any Banach algebra $\mc A$ with the max
quantisation, as in this case, $\wap(\mc A')'$ will also get the max
quantisation.  However, the max quantisation only respects quotients and
not necessarily subspaces.}.

In the case of $L_1(G)$, we can argue as follows, however.  Let $M(G)$
be the Banach algebra of measures on $G$, with convolution product.
Then $M(G) = C_0(G)'$ is a dual Banach algebra, and when we give $M(G)$
the natural operator space structure this induces, we see that
$M(G)$ gets the max quantisation, and is hence a CC Banach algebra.
Then $L_1(G)$ is completely isometrically a subspace of $M(G)$, and
by Theorem~\ref{Main_Thm}, $M(G)$ is completely isometrically a subspace
of $\mc{CB}(E)$ for some reflexive operator space $E$.  The same hence
applies to $L_1(G)$, so by Theorem~\ref{wap_thm}, $\wap(L_1(G)')$ induces
the operator space structure on $L_1(G)$.

Similarly, let $B(G)$ be the Fourier-Stieltjes algebra of $G$, so
that $B(G) = C^*(G)'$, and hence inherits an operator space structure
turning it into a CC dual Banach algebra (see \cite[Section~6]{Runde2}
for example).  Then $A(G)$ is a closed ideal in $B(G)$, and the
operator space structures agree, so we again see that the canonical
map $A(G) \rightarrow \wap(A(G)')'$ is a complete isometry onto its
range.  Hence $\wap(A(G)')$ cannot be too ``badly behaved''.

This paper grew out of an attempt to use the operator space structure
to study weakly almost periodic functionals on $A(G)$.  Our hope was that
using the factorisation definition of weakly compact, we might find a
new definition for operator spaces.  However, as we have seen above,
this is not the case.  Of course, it remains possible that $\wap(A(G)')$
is a well-behaved space, and that simply further work is required.
Alternatively, maybe we need to use the operator space structure on
$A(G)$ in another way.

Two possibilities come to mind.  Let $T:E\rightarrow F$ be a weakly-compact,
completely bounded map.  Then $(T)_n$ is also weakly-compact.  To follow
the analogy with complete boundedness, we would want to attach some value,
corresponding somehow to a ``measure of weak compactness'', to each $(T)_n$,
and define $T$ to be ``completely weakly-compact'' if these values
remained bounded.  Perhaps we could use the factorisation definition of
$T$, and use some invariant of the arising reflexive operator space.

We can define $T:E\rightarrow F$ to be completely bounded if and only if the map
\[ I\otimes T : \mc K(\ell_2) \otimes_{\min} E \rightarrow
\mc K(\ell_2) \otimes_{\min} F \]
is bounded (this approach is taken in \cite{Pisier}).
Obviously $I\otimes T$ is never weakly-compact (even if
$E=F=\mathbb C$, we just get the identity map on $\mc K(\ell_2)$,
which is not a reflexive space).  Can we find some property of
$I\otimes T$ which implies that $T$ is weakly-compact?

\vspace{5ex}

\noindent\emph{Author's Address:}
\parbox[t]{3in}{St. John's College,\\
Oxford,\\
OX1 3JP.}

\bigskip\noindent\emph{Email:} \texttt{matt.daws@cantab.net}

\end{document}